\documentclass[12pt,leqno]{article}
\usepackage{amsmath,amssymb,amscd,latexsym}
\usepackage[all]{xy}
\newcommand{\C}{{\mathbb{C}}}
\newcommand{\F}{{\mathbb{F}}}
\newcommand{\Q}{{\mathbb{Q}}}
\newcommand{\R}{{\mathbb{R}}}
\newcommand{\Z}{{\mathbb{Z}}}
\newcommand{\id}{\mathrm{id}}
\newcommand{\ord}{\mathrm{ord}}
\newcommand{\pr}{\mathrm{pr}}
\newcommand{\res}{\mathrm{res}}
\newcommand{\Aut}{\mathrm{Aut}\,}
\newcommand{\Gal}{\mathrm{Gal}}
\newcommand{\Imm}{\mathrm{Im}\,}
\newcommand{\Ker}{\mathrm{Ker}\,}
\newcommand{\rank}{\mathrm{rank}}
\newcommand{\RRe}{\mathrm{Re}\,}
\newcommand{\sign}{\mathrm{sign}}
\newcommand{\ssp}{\mathrm{span}}
\newcommand{\tors}{\mathrm{tors}}

\newcommand{\Th}{{\mathcal T}}
\newcommand{\eg}{\mathfrak{q}}
\newcommand{\eo}{\mathfrak{o}}
\newcommand{\ep}{\mathfrak{p}}
\newcommand{\oH}{\overline{H}}
\newcommand{\oep}{\overline{\ep}}
\newcommand{\oQ}{\overline{\Q}}
\newcommand{\ox}{\overline{x}}
\newcommand{\osigma}{\overline{\sigma}}
\newcommand{\otau}{\overline{\tau}}
\newcommand{\silo}{\stackrel{\sim}{\longrightarrow}}
\newcommand{\tei}{\, | \,}
\newcommand{\verk}{\mbox{\scriptsize $\,\circ\,$}}
\newcommand{\halb}{\frac{1}{2}}

\newcommand{\te}{\textstyle}

\newtheorem{theorem}{Theorem}[section]

\newtheorem{prop}[theorem]{Proposition}
\newtheorem{cor}[theorem]{Corollary}

\newenvironment{rems}{\bigskip{\bf Remarks}}{}
\newtheorem{conj}[theorem]{Conjecture}

\newtheorem{punkt}[theorem]{$\!\!$}

\newenvironment{proof}{\bigskip{\bf Proof}}{\mbox{}\hfill$\Box$}
\begin{document}
\title{Exchanging the places $p$ and $\infty$ in the Leopoldt conjecture} 
\author{Christopher Deninger}
\date{\ }
\maketitle

\section{Introduction}
\label{sec:1}

It is often reasonable to treat the finite and the infinite places of a number field $k$ on an equal footing.\\
The unit group $E = E_{\infty} (k)$ of $k$ consists of all elements $x$ in $k^*$ with normalized absolute value $\| x \|_v = 1$ for every place $v \nmid \infty$. Fixing a prime $p \neq \infty$ of $\Q$ we can equally well consider the group $E_p (k)$ of all $x$ in $k^*$ with $\| x \|_v = 1$ for all places $v \nmid p$. The elements of $E_p (k)$ are rational $p$-Weil numbers and it is known that $E_p (k)$ is finitely generated. The rank of $E_p (k)$ is known and also a way to construct a subgroup of finite index.

The easier part of the Dirichlet unit theorem asserts that the $(\log \|\;\|_v)$-map 
\begin{equation}
  \label{eq:1}
  E_{\infty} (k) \longrightarrow H_{\infty} = \bigoplus_{v \tei \infty} \R \cdot v
\end{equation}
induces an $\R$-linear injection:
\begin{equation}
  \label{eq:2}
  E_{\infty} (k) \otimes_{\Z} \R \hookrightarrow H_{\infty} \; .
\end{equation}
For $E_p (k)$ there is a corresponding map constructed by Gross \cite{G}:
\begin{equation}
  \label{eq:3}
  E_p (k) \longrightarrow H_p = \bigoplus_{v\tei p} \Q_p \cdot v \; .
\end{equation}
In this case the injectivity of the analogue of (\ref{eq:2})
\begin{equation}
  \label{eq:4}
  E_p (k) \otimes_{\Z} \Q_p \hookrightarrow H_p
\end{equation}
is only a conjecture. Gross proved it under certain conditions using Brumer's $p$-adic analogue of Baker's result on linear independence of logarithms of algebraic numbers. See \cite{G} Cor. 2.14 or section two below.

Leopoldt's conjecture is concerned with the $p$-adic logarithm
\begin{equation}
  \label{eq:5}
  E_{\infty} (k) \longrightarrow k_p = k \otimes \Q_p \; .
\end{equation}
Again the assertion is that its $\Q_p$-linear extension
\begin{equation}
  \label{eq:6}
  E_{\infty} (k) \otimes_{\Z} \Q_p \longrightarrow k_p
\end{equation}
is injective. Again using Brumer's result this can be shown for abelian extensions $k / \Q$. \\
In this note we will first review the preceeding facts in more detail. Then we consider the multivalued $\log$-map
\begin{equation}
  \label{eq:7}
  E_p (k) \longrightarrow k_{\infty} = k \otimes \R 
\end{equation}
which is a counterpart to the Leopoldt map (\ref{eq:5}) with the places $p$ and $\infty$ exchanged.\\
As it turns out the analogue of the Leopoldt conjecture predicts that certain  vectors should be $\R$-linear independent whose components are arguments of conjugates of Weil numbers in $k$. We give the precise formulation in \ref{t31}. It is more involved than for the maps (\ref{eq:2}), (\ref{eq:4}) and (\ref{eq:6}) since the complex logarithm is multivalued. Using Baker's result on linear forms in logarithms we prove part of our conjecture in certain abelian situations.

The multivalued map (\ref{eq:7}) gives rise to a single-valued map of $E_p (k)$ into a real torus. In section \ref{sec:4} we describe the closure of its image.

The Dirichlet regulator constructed via (\ref{eq:2}) is related to the zeta function of $k$ both at $s = 0$ and $s = 1$. The Gross regulator based on (\ref{eq:4}) is related to the $p$-adic zeta function of $k$ at $s = 0$. The Leopoldt regulator based on (\ref{eq:6}) is related to its behaviour at $s = 1$. It is an open question whether the regulators one can construct using (\ref{eq:7}) are connected to some kind of zeta function. 

The relation of Weil numbers with motives over finite fields is explained in \cite{M}. In \cite{GO}, \cite{O} and \cite{K}, Weil numbers are studied from the point of view of algebraic and analytic number theory. 

It is a pleasure for me to thank Niko Naumann and Kay Wingberg for helpful comments and the CRM in Montreal for its support. 

\section{Review of some material from algebraic number theory}
\label{sec:2}

For a number field $k / \Q$ and a prime number $p$ we review some known facts about the group $ E_p (k)$ defined in the introduction. Our basic reference is \cite{M}. \\
An algebraic number $\lambda$ is a $q = p^r$-Weil number of weight $w$ if it has the following two properties:

i) For every embedding $\sigma : \Q [\lambda] \hookrightarrow \C$ we have $|\sigma (\lambda)| = q^{w/2}$.\\
ii) There is some $n \in \Z$ such that $p^n \lambda$ is an algebraic integer. 

The set $W_w (q)$ of Weil numbers of weight $w$ in an algebraic closure $\oQ$ of $\Q$ is stable under the action of $G_{\Q} = \Gal (\oQ / \Q)$. The field generated by a set of $q$-Weil numbers is either totally real or a $CM$ field \cite{GO} or \cite{M} Def. 2.5 ff.

If $k$ contains a $CM$ field, let $k_{CM}$ denote the maximal $CM$ field contained in $k$. The following is clear:

\begin{prop}
  \label{t21}
1) For any embedding $k \subset \oQ$ we have:
\[
E_p (k) = W_0 (p) \cap k^* \; .
\]
2) If $k$ contains a $CM$ field, then
\[
E_p (k) = E_p (k_{CM}) \; .
\]
If not then $E_p (k) = \{ \pm 1 \}$.
\end{prop}

In the study of the unit group $E = E_{\infty} (k)$ one considers the regulator map obtained by composition
\begin{equation}
  \label{eq:8}
  \alpha_{\infty} : E_{\infty} (k) \xrightarrow{\Delta} \bigoplus_{v\tei \infty} k^*_v \xrightarrow{(\| \; \|)_{v \tei \infty}} \bigoplus_{v\tei \infty} \| k^*_v \|_v \xrightarrow{\bigoplus_v \log} \bigoplus_{v\tei \infty} \R \cdot v \; .
\end{equation}
Here $\Delta$ is the diagonal inclusion, $k_v$ is the $v$-completion and $\| k^*_v \|_v$ is the value group $\R^*_+$ of the absolute value $\| \; \|_v$. An analogous map is used for the study of $E_p (k)$ in \cite{M} Prop. 2.27. 
It is the composition
\begin{equation}
  \label{eq:9}
  \alpha_{(p)} : E_p (k) \xrightarrow{\Delta} \bigoplus_{v \tei p} k^*_v \xrightarrow{(\|\; \|_v)_{v \tei p}} \bigoplus_{v \tei p} \| k^*_v \|_v \xrightarrow{\bigoplus_v \log_{(p)}} \bigoplus_{v \tei p} \Z \cdot v \; .
\end{equation}
Here $\log_{(p)} (r) = \log r / \log p$ for $r \in \| k^*_v \|_v = Nv^{\Z}$. Thus $\alpha_{(p)}$ is given by
\begin{equation}
  \label{eq:10}
  \alpha_{(p)} (x) = - \sum_{v\tei p} f (v \tei p) \ord_v (x) \cdot v \; .
\end{equation}
By \ref{t21} we may assume that $k$ is a $CM$ field. We then denote the places dividing $p$ by $\ep$ instead of $v$. Let $c$ denote complex conjugation on $k$. Then we have for $x \in E_p (k)$ that
\[
\ord_{\ep} (x) + \ord_{\ep^c} (x) = \ord_{\ep} (x) + \ord_{\ep} (x^c) = \ord_{\ep} (xx^c) = \ord_{\ep} (1) = 0 \; .
\]
Thus $\alpha_{(p)}$ takes values in
\[
\Big( \bigoplus_{\ep \tei p} \Z \cdot \ep \Big)^- = \Big\{ \sum_{\ep \tei p} \nu_{\ep} \cdot \ep \, \Big| \, \nu_{\ep} + \nu_{\ep^c} = 0 \quad \mbox{for all} \; \ep \Big\} \; .
\]
Besides the map
\begin{equation}
  \label{eq:11}
  \alpha_{(p)} : E_p (k) \longrightarrow \Big( \bigoplus_{\ep \tei p} \Z \cdot \ep \Big)^-
\end{equation}
there are maps in the other direction. Let $\Th$ be the set of primes $\ep \tei p$ in $k$ with $\ep \neq \ep^c$. Choose an integer $M \ge 1$ which is divisible by all inertia indeces $f (\ep \tei p)$ for $\ep$ in $\Th$ and such that $\ep^{M f(\ep \tei p)^{-1}}$ is a principal ideal for these $\ep$. We may then choose $x_{\ep}$ in $\eo_k$ for $\ep$ in $\Th$ such that we have $\ep^{Mf (\ep \tei p)^{-1}} = (x_{\ep})$ and $x^c_{\ep} = x_{\ep^c}$. Define a map
\begin{equation}
  \label{eq:12}
  \pi_M : \Big( \bigoplus_{\ep \tei p} \Z \cdot \ep \Big)^- \longrightarrow E_p (k)
\end{equation}
by the formula
\[
\pi_M \Big( \sum_{\ep \tei p} \nu_{\ep} \cdot \ep \Big) = \prod_{\ep\in \Th} x^{\nu_{\ep}}_{\ep} \; .
\]
Essentially, the following facts can be found in \cite{S} Ch. II or in \cite{M}. For convenience we give the short proof. 

\begin{prop}
  \label{t22}
Let $k$ be a $CM$ field and let $S$ be a system of representatives of the primes in $\Th$ under complex conjugation. The following assertions hold:\\
1) $\alpha_{(p)} \verk \pi_M = - M \cdot \id$ and $(\pi_M \verk \alpha_{(p)}) (x) \equiv x^{-M} \mod \mu (k)$.\\
2) The abelian group $E_p (k)$ is finitely generated with
\begin{eqnarray*}
  \tors E_p (k) & = & \mu (k) = \Ker \alpha_{(p)} \quad \mbox{and} \\
\rank E_p (k) & = & \halb |\Th| \le \halb (k : \Q) \; .
\end{eqnarray*}
3) Rationally $\alpha_{(p)}$ and $\pi_M$ are isomorphisms.\\
4) The elements $\xi_{\ep} = x^c_{\ep} x^{-1}_{\ep} \in E_p (k)$ for $\ep \in S$ form a $\Q$-basis of $E_p (k) \otimes \Q$.
\end{prop}

\begin{proof}
  1) The first relation follows from a straightforward computation. As for the second, set
\[
y = x^M \pi_M (\alpha_{(p)} (x)) = x^M \prod_{\ep \in \Th} x^{-f (\ep \tei p) \ord_{\ep} (x)}_{\ep} \; .
\]
Then $(y) = (x)^M \prod_{\ep \tei p} \ep^{-M \ord_{\ep} (x)} = (1)$. Hence $y$ is a unit.

Since $y$ lies in $E_p (k)$ it follows that all absolute values of $y$ are equal to $1$. Thus $y$ is a root of unity.\\
2, 3) Since $\alpha_{(p)}$ embeds $E_p (k) / \mu (k)$ injectively into a finitely generated abelian group, $E_p (k)$ is itself finitely generated. By 1) the map $\alpha_{(p)}$ is an isomorphism rationally, hence
\[
\rank \, E_p (k) = \dim \Big( \bigoplus_{\ep \tei p} \Q \cdot \ep \Big)^- = \halb |\Th| \; .
\]
4) The group $\left( \bigoplus_{\ep \tei p} \Z \cdot \ep \right)^-$ has basis the elements $\ep^c - \ep$ for $\ep \in S$. They are mapped to $x_{\ep^c} x^{-1}_{\ep} = x^c_{\ep} x^{-1}_{\ep}$ under $\pi_M$. Now 3) implies the assertion.
\end{proof}

Apart from (\ref{eq:9}) there is a deeper $p$-analogue of the Dirichlet regulator map $\alpha_{\infty}$ in (\ref{eq:8}). It was introduced by Gross in order to formulate a $p$-adic analogue of the Stark conjectures \cite{G}.

Gross first defines ``absolute values'' with values in $\Z_p$'' for every valuation $v$ of $k$.
Namely, for $x$ in $k^*_v$ set
\[
\| x \|_{v,p} = \left\{ 
  \begin{array}{ll}
1 & v \; \mbox{complex} \\
\sign \, x & v \; \mbox{real} \\
(Nv)^{-\ord_v (x)} & v \; \mbox{finite} \; v \nmid p \\
(Nv)^{-\ord_v (x)} N_{k_v / \Q_p} (x) & v \tei p \; .
  \end{array} \right.
\]
As for the usual absolute values the product formula holds:
\[
\prod_v \| x \|_{v,p} = 1 \quad \mbox{for} \; x \in k^* \; .
\]
Let $\log_p : \Z^*_p \to \Q_p$ be the $p$-adic logarithm. Then Gross' $p$-adic analogue of $\alpha_{\infty}$ is the composition:
\begin{equation}
  \label{eq:13}
  \alpha_p : E_p (k) \xrightarrow{\Delta} \bigoplus_{v \tei p} k^*_v \xrightarrow{(\|\; \|_{v,p})_{v \tei p}} \bigoplus_{v\tei p} \| k^*_v \|_{v,p} \xrightarrow{\bigoplus_v \log_p} \bigoplus_{v \tei p} \Q_p \cdot v \; .
\end{equation}
Its kernel consists of the roots of unity in $k$ and Gross conjectures that tensored with $\Q_p$ the map $\alpha_p$ remains injective, \cite{G} conjecture 1.15. 

In \cite{G} corollary 2.14 this is proved in the following case: $k$ contains a $CM$ field and the set of places $\ep^+ \tei p$ in $k^+_{CM}$ which split in $k_{CM}$ is permuted transitively by an abelian subgroup of $\Aut k^+_{CM}$.

We now turn to the Leopoldt conjecture whose analogue for $p$ and $\infty$ exchanged will be studied in the next section.

For every place $v \tei p$ of $k$ let $\log_v : k^*_v \to k_v$ be the $p$-adic logarithm normalized by $\log_v p = 1$. Consider the composition:
\begin{equation}
  \label{eq:14}
  \alpha_{\infty , p} : E_{\infty} (k) \xrightarrow{\Delta} \bigoplus_{v\tei p} k^*_v \xrightarrow{\bigoplus_v \log_v} \bigoplus_{v \tei p} k_v = k \otimes \Q_p \; .
\end{equation}
Its kernel consists of roots of unity, so that
\[
\alpha_{\infty , p} \otimes 1 : E_{\infty} (k) \otimes \Q \longrightarrow k \otimes \Q_p
\]
is injective. The Leopoldt conjecture asserts that the map 
\[
\alpha_{\infty , p} \otimes 1 : E_{\infty} (k) \otimes \Q_p \longrightarrow k \otimes \Q_p 
\]
is injective. It is known for abelian extensions $k / \Q$ for example \cite{W}.
\section*{Appendix}
Let $\lambda$ be a $q$-Weil number of weight $w$ in $\oQ \subset \C$. Let $k \subset \C$ be a $CM$ field containing $\lambda$. Then $\lambda$ is determined up to a root of unity by the valuations $\|\lambda \|_{\ep}$ for $\ep \tei p$ in $k$. 

If $\alpha \in \C$ is such that $q^{\alpha} = \lambda$, then $\RRe \alpha = \frac{w}{2}$. The imaginary part of $\alpha$ is determined by $\lambda$ up to an element of $\frac{2\pi}{\log q} \Z$. Hence the values $\| \lambda \|_{\ep}$ for $\ep \tei p$ determine $\Imm \alpha \mod \frac{2\pi}{\log q} \Q$. The material of the first section allows one to write down a formula for this dependence. For a variety $X / \F_q$ this gives a way to reflect $\ep$-adic divisibility properies of the zeroes of the zeta function $Z_X (T)$ in the imaginary parts of the zeroes of the Hasse Weil zeta function $\zeta_X (s)$. Recall here that
\[
\zeta_X (s) = \prod_{x \in |X|} (1 - N x^{-s})^{-1} \quad \mbox{for} \; s > \dim X \quad \mbox{in} \; \C
\]
and that the two zeta functions are related by the formula:
\[
\zeta_X (s) = Z_X (q^{-s}) \; .
\]
In particular, I wonder whether a statement such as ``the Newton polygon lies above the Hodge polygon'' might be reformulated in terms of $\zeta_X (s)$ and its zeroes such that it makes sense for $X / \Z$.

However that may be, the desired formula for $\Imm \alpha$ in $\R / \frac{2\pi}{\log q} \Q$ is the following:
\begin{eqnarray}
  \label{eq:15}
  \widetilde{\Imm} \alpha & = & \sum_{\ep \in \Th} \log_{(p)} \| \lambda \|_{\ep} \; \widetilde{\arg}_q (x^{-1/M}_{\ep}) \\
& = & \sum_{\ep \in \Th} f (\ep \tei p) \ord_{\ep} (\lambda) \widetilde{\arg}_q (x^{1/M}_{\ep}) \; . \nonumber
\end{eqnarray}
Here $\arg_q = (\log q)^{-1} \arg$ and for a real number $t$ we let $\tilde{t}$ be its image in $\R / \frac{2\pi}{\log q} \Q$. The integer $M$ is chosen as in section 2 before proposition \ref{t22}.

\begin{proof}
  The number $x = \lambda / \lambda^c$ is in $E_p (k)$. By proposition \ref{t22} we have
  \begin{eqnarray}
    2 \, \widetilde{\arg}_q \lambda & = & \widetilde{\arg}_q x = -M^{-1} \widetilde{\arg} (\pi_M \alpha_{(p)} (x)) \nonumber \\
& = & - M^{-1} \widetilde{\arg}_q \pi_M \left( \sum_{\ep \in \Th} \log_{(p)} \| \lambda / \lambda^c \|_{\ep} \cdot \ep \right) \nonumber \\
& = & - M^{-1} \sum_{\ep \in \Th} \log_{(p)} \| \lambda / \lambda^c \|_{\ep} \cdot \widetilde{\arg}_q x_{\ep} \; . \label{eq:16}
  \end{eqnarray}
Since $x_{\oep} = \ox_{\ep}$ we find
\[
\widetilde{\arg}_q \lambda = - M^{-1} \sum_{\ep \in \Th} \log_{(p)} \| \lambda \|_{\ep} \; \widetilde{\arg}_q x_{\ep} \; .
\]
Hence (\ref{eq:15}) follows.
\end{proof}

Another way to write $\widetilde{\Imm} \alpha$ is the following: Let $S$ be as in proposition \ref{t22}. Then starting from (\ref{eq:16}) we find
\[   
\widetilde{\Imm} \alpha = \sum_{\ep \in S} \log_{(p)} (\|\lambda \|_{\ep} / \| \lambda \|_{\ep^c}) \; \widetilde{\arg}_q (\xi_{\ep}) \quad \mbox{in} \; \R / {\te \frac{2\pi}{\log q}} \Q
\]
where $\xi_{\ep} = (x^c_{\ep} x^{-1}_{\ep})^{1/M}$ in $\C^* / \mu (\C)$.

In the $\Q$-vector space $\R / \frac{2\pi}{\log q} \Q$ the vectors $\widetilde{\arg}_q (\xi_{\ep})$ for $\ep$ in $S$ are $\Q$-linearly independent by Prop. \ref{t22}, 4). Hence we see how $\widetilde{\Imm} \alpha$ determines the quotients $\| \lambda \|_{\ep} / \| \lambda \|_{\ep^c}$ for all $\ep \in S$. These in turn determine $\| \lambda \|_{\ep}$ for all $\ep \tei p$ since $\lambda$ was a $q$-Weil number of weight $w$. 

\section{Changing places in the Leopoldt conjecture}
\label{sec:3}

In this section we deal with an analogue of the Leopoldt regulator map (\ref{eq:14}) for $E_p (k)$ instead of $E_{\infty} (k)$. We may therefore assume that $k$ is a $CM$ field. Ideally we would consider the composition:
\begin{equation}
  \label{eq:17}
  E_p (k) \xrightarrow{\Delta} \bigoplus_{v \tei \infty} k^*_v \xrightarrow{\bigoplus_v \log_v} \bigoplus_{v \tei \infty} k_v = k \otimes \R \; .
\end{equation}
Here $\log_v$ is the complex logarithm on $k^*_v \cong \C^*$. Unfortunately the complex logarithm is a multivalued map, so (\ref{eq:17}) has to be interpreted suitably.

Note that for $x$ in $E_p (k)$ we have
\[
2 \RRe \log_v x = \log \| x \|_v = \log 1 = 0 \; .
\]
We set
\[
\R (1)_v = \{ z \in k_v \tei \RRe z = 0 \} \; .
\]
More generally, for any subgroup $\Lambda \subset \R$ let us write $\Lambda (1)_v$ for the image of $2 \pi i \Lambda \subset \C$ in $k_v$ under any of the two continuous isomorphisms $\C \silo k_v$. Then we interpret (\ref{eq:17}) as the injective map:
\begin{equation}
  \label{eq:18}
  \alpha_{p,\infty} : E_p (k) \xrightarrow{\Delta} \bigoplus_{v \tei \infty} k^*_{v,1} \xrightarrow{\bigoplus_v \log_v} T_k = \bigoplus_{v \tei \infty} \R (1)_v / \Z (1)_v \; .
\end{equation}
Here $k^*_{v,1}$ is the subgroup of $k^*_v$ of elements with norm one.

Since $T_k$ is not an $\R$-vector space we cannot state a ``Leopoldt conjecture'' in quite the same way as before. The strongest formulation seems to be the following.

\begin{conj}
  \label{t31}
Let $\xi_1 , \ldots , \xi_m$ be a basis of $E_p (k) \otimes \Q$ and consider the diagram
\[
E_p (k) \otimes \Q \xrightarrow{\alpha_{p,\infty} \otimes 1} \bigoplus_{v\tei \infty} \R (1)_v / \Q (1)_v \xleftarrow{\pr} V_k = \bigoplus_{v \tei \infty} \R (1)_v \; .
\]
For any $1 \le i \le m$ choose a preimage $\eta_i$ under $\pr$ of $(\alpha_{p,\infty} \otimes 1) (\xi_i)$. Then the vectors $\eta_i$ are $\R$-linearly independent in $V_k$.
\end{conj}

Note here that
\[
\rank \; E_p (k) = |S| \le  \halb (k : \Q) = \dim V_k \; .
\]
We have equality if and only if $p$ splits completely in $k$. In this case one may form regulator determinants. Observe also that we get a weaker conjecture if we work integrally instead of tensoring with $\Q$.

\begin{prop}
  \label{t32}
It is sufficient to verify the conjecture for any fixed basis $\xi_1 , \ldots , \xi_m$ (but all possible liftings $\eta_i$).
\end{prop}

The proof is straightforward.

Choose a representative $\sigma_v : k \hookrightarrow \C$ for every complex place $v$ of $k$. Let $S$ be as in proposition \ref{t22}. Then according to propositions \ref{t22}, 4) and \ref{t32} our conjecture is equivalent to the following:

\begin{conj}
  \label{t33} The vectors $(\arg (\overline{\sigma_v (x_{\ep})} \sigma_v (x_{\ep})^{-1}))_{v \tei \infty}$ for $\ep \in S$ are $\R$-linearly independent in $\R^{r_2}$ for any choices of argument $\mod 2 \pi \Q$.
\end{conj}

{\bf Trivial example} Assume that there are exactly two primes dividing $p$ in $k$. They are $\ep $ and $\ep^c \neq \ep$ for some $\ep \tei p$. Take $S = \{ \ep \}$. The conjecture asserts that for some $v \tei \infty$
\[
\arg (\overline{\sigma_v (x_{\ep})} \sigma_v (x_{\ep})^{-1}) \notin 2 \pi \Q \; .
\]
In other words $\overline{\sigma_v (x_{\ep})} \sigma_v (x_{\ep})^{-1}$ should not be a root of unity. This is clear, since otherwise we would have $(\ep^c)^N = \ep^N$ for some $N \ge 1$, a contradiction.

Next we introduce special bases of $E_p (k)$, at least in favourable circumstances.

\begin{prop}
  \label{t34}
Let $k$ be a $CM$ field and consider as before the set $\Th$ of primes $\ep$ dividing $p$ in $k$ such that $\ep^c \neq \ep$. Assume that $\Aut k$ contains a subgroup $H$ with $c \notin H$ such that $G = \langle H , c \rangle$ permutes the primes in $\Th$ transitively with isotropy subgroups contained in $H$.\\
Then there exists an element $\xi$ in $E_p (k)$ such that the set of conjugates $\xi^{\sigma}$ for $\sigma$ in $H$ is a basis of $E_p (k) \otimes \Q$.
\end{prop}

\begin{rems}
  {\bf 1} The conjugates $\xi^{\sigma}$ do not need to be pairwise different.\\
{\bf 2} As $c$ is central we have $G \cong H \times \langle c \rangle$ and $H$ is normal in $G$. Hence the condition on isotropy subgroups has to be checked for one prime in $\Th$ only. It is trivially satisfied if $p$ is completely decomposed in $k$. In that case $k / \Q$ must be Galois and the conjugates of $\xi$ under $H$ are pairwise distinct.
\end{rems}

\begin{proof}
  We may assume that $\Th$ is non-empty. Fix a prime $\ep_0$ in $\Th$, let $k'$ be the fixed field of $G$ and set $\ep'_0 = \ep_0 \cap k'$. Then we claim:
  \begin{equation}
    \label{eq:19}
    \Th = \; \mbox{set of primes} \; \ep \; \mbox{in} \; k \; \mbox{dividing} \; \ep'_0 \; .
  \end{equation}
Indeed, if $\ep \neq \ep^c$ divides $p$, then by assumption on $G$ there is some $\sigma \in G$ with $\ep = \ep^{\sigma}_0$. Hence $\ep \cap k' = (\ep_0 \cap k')^{\sigma} = \ep'_0$. On the other hand, the prime divisors $\ep$ of $\ep'_0$ are all conjugate under $G$ since the extension $k / k'$ is Galois with group $G$. For any $\ep \tei \ep'_0$ there is thus some $\sigma$ in $G$ with $\ep = \ep^{\sigma}_0$. Hence $\ep^c = \ep^{\sigma c}_0 = \ep^{c\sigma}_0 \neq \ep^{\sigma}_0 = \ep$ since $\ep^c_0 \neq \ep_0$. 

By proposition \ref{t22} the map $\alpha_{(p)}$ in (\ref{eq:11}) induces an isomorphism
\begin{equation}
  \label{eq:20}
  \alpha_{(p)} : E_p (k) \otimes \Q \silo \Big( \bigoplus_{\ep \tei p} \Q \cdot \ep \Big)^- = \Q [\Th]^-
\end{equation}
where $(\;)^-$ denotes the $-1$ eigenspace of the action of complex conjugation. The map $\alpha_{(p)}$ is $G$-equivariant under the natural left $G$-actions on both sides.

Now let $G_{\ep_0}$ be the decomposition group of $\ep_0$ in the extension $k / k'$. According to (\ref{eq:19}) we have a bijection
\begin{equation}
  \label{eq:21}
  G / G_{\ep_0} \silo \Th \; , \; \sigma G_{\ep_0} \longmapsto \sigma (\ep_0) \; .
\end{equation}
Together with (\ref{eq:20}) we obtain a $G$-equivariant isomorphism:
\begin{equation}
  \label{eq:22}
  E_p (k) \otimes \Q \silo \Q [G / G_{\ep_0}]^- \; .
\end{equation}
By assumption $G_{\ep_0} \subset H$. Hence we get an embedding
\[
H / G_{\ep_0} \hookrightarrow G / G_{\ep_0} \cong \langle c \rangle \times H / G_{\ep_0} \; .
\]
The natural restriction map
\[
\res : \Q [G / G_{\ep_0}] \longrightarrow \Q [H / G_{\ep_0}]
\]
induces an $H$-equivariant isomorphism
\[
\res : \Q [G / G_{\ep_0}]^- \silo \Q [H / G_{\ep_0}] \; .
\]
Combined with (\ref{eq:22}) we get an $H$-equivariant isomorphism:
\begin{equation}
  \label{eq:23}
  E_p (k) \otimes \Q \silo \Q [H / G_{\ep_0}] \; .
\end{equation}
Hence any element $\xi$ in $E_p (k)$ corresponding to a multiple of $\delta_e = 1 \cdot e G_{\ep_0}$ has the required property. 

More {\it explicitely} we can describe such a $\xi$ as follows. Let $k_{\ep_0}$ be the fixed field of $G_{\ep_0}$. In the extension $k / k_{\ep_0}$ only one prime lies above $\eg_0 = \ep_0 \cap k_{\ep_0}$, namely $\ep_0$. Thus we have $\eg_0 \eo_k = \ep^e_0$ where $e \ge 1$. Choose $y \in k_{\ep_0}$ such that $(y) = \eg^h_0$ in $k_{\ep_0}$ for some $h \ge 1$. It follows that $(y) = \ep^{eh}_0$ in $k$. Choose a multiple $M \ge 1$ of $efh$ where $f$ is the common inertia index over $p$ of the primes in $\Th$ and set $x = y^{H / efh} \in k_{\ep_0}$. Then we have $(x) = \ep^{Hf^{-1}}_0$ in $k$. The $H / G_{\ep_0}$ conjugates of $\xi = x^c x^{-1} \in E_p (k_{\ep_0}) \subset E_p (k)$ give a basis of $E_p (k) \otimes \Q$.
\end{proof}

\begin{punkt}
  \label{t36} \rm
We now prove special cases of conjecture \ref{t31}. Let $k$ be a $CM$ field and assume that we are given an element $\xi$ in $E_p (k)$ and an abelian subgroup $H$ of $\Aut k$ with $c \notin H$ such that the set of conjugates $\xi^{\sigma}$ for $\sigma \in H$ is a basis of $E_p (k) \otimes \Q$.

Let $H_{\xi}$ be the group of $\sigma \in H$ fixing $\xi$ and let $\oH$ be the abelian quotient group $\oH = H / H_{\xi}$. Fix an embedding $k \subset \C$.

For every $\osigma \in \oH$ choose a preimage denoted by $\sigma_v$ in $H$ and hence in $\Aut k$. Via the embedding $k \subset \C$ we will also view $\sigma_v$ as an embeding $\sigma_v : k \hookrightarrow \C$. Let $v$ be the infinite place corresponding to $\sigma_v$. Then $\sigma_v$ induces isomorphisms $\sigma_v : k_v \silo \C$ and $\sigma_v : \R (1)_v \silo \R (1)$. Let $S_{\infty}$ be the set of infinite places so obtained. The bijection $\oH \silo \{ \xi^{\sigma} \tei \sigma \in H \} , \otau \mapsto \otau (\xi)$ shows that 
\[
|S_{\infty}| = |\oH| = \rank E_p (k) \; .
\]
Next, let $s$ be a set theoretical splitting of the projection $\R (1) \to \R (1) / \Q (1)$. It induces splittings $s_v$ of $\R (1)_v \to \R (1)_v / \Q (1)_v$ via conjugation, $s_v = \sigma^{-1}_v s \sigma_v$.

For each basis vector $\otau (\xi)$ of $E_p (k) \otimes \Q$ with $\otau \in \oH$ define as follows a lift $\eta_{\otau}$ to $V_k$ of the vector
\[
(\alpha_{p,\infty} \otimes 1) (\otau (\xi)) = (\log_v (\otau (\xi)_v))_{v\tei \infty} \; .
\]
The $v$-components with $v$ in $S_{\infty}$ are lifted via $s_v$. The other components are lifted arbitrarily. For $v$ in $S_{\infty}$ we therefore have:
\[
\eta_{\otau , v} = s_v \log_v (\otau (\xi)_v) \; .
\]
\end{punkt}

\begin{theorem}
  \label{t37}
In the situation of \ref{t36} the vectors $\eta_{\otau}$ for $\otau \in \oH$ are $\R$-linearly independent in $V_k$.
\end{theorem}

\begin{proof}
  It suffices to show that
\[
\delta = |\det (\sigma_v s_v \log_v (\otau (\xi)_v))_{\otau \in \oH , v \in S_{\infty}}|
\]
is non-zero. Note that the determinant itself is well defined only up to sign.

We have
\begin{eqnarray*}
  \sigma_v s_v \log_v (\otau (\xi)_v) & = & s \sigma_v \log_v (\otau (\xi)_v) = s \log \sigma_v \otau (\xi) \\
& = & s \log \osigma_v \otau (\xi) \; .
\end{eqnarray*}
We may order $S_{\infty}$ and $\oH$ in such a way that we get:
\[
\delta = |\det (s \log \osigma^{-1} \otau (\xi))_{\osigma , \otau \in \oH}|\; .
\]
The well known formula for the group determinant \cite{W} Lemma 5.26 gives:
\[
\delta = \prod_{\chi \in \oH^*} \Big| \sum_{\osigma \in \oH} \chi (\osigma) s \log \osigma (\xi) \Big| \; .
\]
We may choose $N \ge 1$ so that
\[
N s\log \osigma (\xi) = \log_{\osigma} (\osigma (\xi)^N) \quad \mbox{for all} \; \osigma \in \oH \; .
\]
Here $\log_{\osigma}$ is a suitable branch of the complex logarithm depending on $\osigma$. It suffices therefore to show that
\[
\sum_{\osigma \in \oH} \chi (\osigma) \log_{\osigma} (\osigma (\xi)^N) \neq 0 \quad \mbox{for every character} \; \chi \; .
\]
This is an immediate consequence of Baker's theorem \cite{B} Theorem 2.1, since the numbers $\log_{\osigma} (\osigma (\xi)^N)$ are $\Q$-linear independent as one sees upon exponentiating.
\end{proof}

\begin{cor}
  \label{t38}
Let $k$ be a $CM$ field satisfying the assumptions of proposition \ref{t34} with $H$ abelian and choose any $\xi$ constructed there. We may then apply the construction of \ref{t36} to get liftings $\eta_{\otau}$ of $(\alpha_{p,\infty} \otimes 1) (\otau (\xi))$ for $\otau \in \oH$ which are linearly independent in $V_k$. 
\end{cor}

This lends some credibility to the conjecture. In general, even in the abelian case it seems beyond the power of Baker's theorem. For example assume we are in the situation of proposition \ref{t34} with $H = \{ 1 , \sigma \}$ of order two and exactly four different primes of $k$ dividing $p$. Then $\xi , \xi^{\sigma}$ form a basis of $E_p (k) \otimes \Q$. Baker's theorem allows us to show that for any two branches $\log^{(1)}$ and $\log^{(2)}$ of the complex logarithm, the determinant
\[
\begin{vmatrix}
  \log^{(1)} \xi & \log^{(2)} \xi^{\sigma} \\
\log^{(2)} \xi^{\sigma} & \log^{(1)} \xi
\end{vmatrix}
\]
is non-zero. However if we choose the branches of the logarithm less coherently we are no longer reduced to {\it linear} forms in logarithms. For example, I cannot prove the following determinant to be non-zero:
\[
\begin{vmatrix}
  \log^{(1)} \xi & \log^{(2)} \xi^{\sigma} \\
\log^{(2)} \xi^{\sigma} & \log^{(3)} \xi
\end{vmatrix}
\]
where $\log^{(3)} = \log^{(1)} + 2 \pi i$.

Let us finally discuss a case where $k / \Q$ is abelian but the requirement in proposition \ref{t34} that $\langle c \rangle$ has a complement is not met. Take $k = \Q (\zeta_5)$. Its Galois group is cyclic of order $4$. Let $\sigma_{\nu}$ be the automorphism corresponding to $\nu$ in $(\Z / 5)^*$.

We are interested in the prime $p = 11$. It decomposes into the primes \\
$\ep_{\nu} = (1 + 2 \zeta^{\nu}_5)$ for $\nu \in (\Z / 5)^*$. Thus $\Th = \{ \ep_{\nu} \tei \nu \in (\Z / 5)^* \}$. Let us choose $S = \{ \ep_1 , \ep_2 \}$ for example. We may take $M = 1$ in the construction before proposition \ref{t22} and
\[
x_{\ep_1} = 1 + 2 \zeta_5 \quad \mbox{and} \quad x_{\ep_2} = 1 + 2 \zeta^2_5 \; .
\]
Then
\[
\xi_{\ep_1} = x^c_{\ep_1} x^{-1}_{\ep_1} \quad \mbox{and} \quad \xi_{\ep_2} = x^c_{\ep_2} x^{-1}_{\ep_2} = \xi^{\sigma_2}_{\ep_1}
\]
form a basis of $E_{11} (k) \otimes \Q$. Now let us use the same branch $\log$ of logarithm for the liftings in conjecture \ref{t31}. Since $\sigma^2_2 = \sigma_{-1}$ instead of $\sigma^2_2 = 1$ we do not end up with a group determinant. Instead we get a determinant like the following
\[
\begin{vmatrix}
  \log \xi_{\ep_1} & \log \xi^{\sigma_2}_{\ep_1} \\
\log \xi^{\sigma_2}_{\ep_1} & \log \xi^{\sigma^2_2}_{\ep_1}
\end{vmatrix} = 
\begin{vmatrix}
  \log \xi_{\ep_1} & \log \xi_{\ep_2} \\
\log \xi_{\ep_2} & \log \xi^{-1}_{\ep_1}
\end{vmatrix} \; .
\]
For an unfortunate choice of $\log$ we have $\log \xi^{-1}_{\ep_1} = - \log \xi_{\ep_1} \pm 2 \pi i$ so that the determinant does not factor into a product of linear forms in logarithms. We have therefore not persued the case where $\langle c \rangle$ does not have a complement $H$.

\section{On the closure of $\alpha_{p,\infty} (E_p (k))$ in $T_k$}
\label{sec:4}

The Leopoldt conjecture can be stated as an assertion about the closure of the global units embedded into the local units at $p$. Regarding the map $\alpha_{p,\infty}$ in (\ref{eq:18}) one may ask about the closure of the subgroup $\alpha_{p,\infty} (E_p (k))$ in the real $\halb (k : \Q)$-dimensional torus $T_k$. We have
\[
\overline{\alpha_{p,\infty} (E_p (k))} = \alpha_{p,\infty} (E_p (k) \cap \mu (k)) \times \overline{\alpha_{p,\infty} (E_p (k))}^0
\]
where $^0$ denotes the connected component of zero. In order to determine $ \overline{\alpha_{p,\infty} (E_p (k))}^0$ it is equivalent but more convenient to describe the quotient
\[
C = T_k /   \overline{\alpha_{p,\infty} (E_p (k))}^0 \; .
\]
For the Pontrjagin dual we have:
\begin{eqnarray*}
  \hat{C} & = & \Ker \left( \hat{T}_k \longrightarrow (\overline{\alpha_{p,\infty} (E_p (k))}^0)^{\wedge} \right) \\
& = & \Ker \left( \hat{T}_k \longrightarrow \alpha_{p,\infty} (E_p (k))^{\wedge} \otimes \Q \right) \; .
\end{eqnarray*}
For simplicity assume that $k / \Q$ is a Galois $CM$ extension with group $G$. Fix an embedding $k \subset \C$ and for every $v \tei \infty$ choose some $\sigma_v \in G$ which gives rise to the place $v$ after composing with the inclusion $k \hookrightarrow \C$.

There is an isomorphism
\[
\bigoplus_{v \tei \infty} \Z \silo \hat{T}_k \; , \; \mu \longmapsto \chi_{\mu}
\]
with $\chi_{\mu} (t) = \prod_{v \tei \infty} \exp \mu_v (\sigma_v t_v)$.

Here $\sigma_v$ is viewed as the isomorphism $\sigma_v : \R (1)_v / \Z (1)_v \to \R (1) / \Z (1)$ induced by $\sigma_v : k_v \silo \C$.

Since
\[
\sigma_v \log_v \xi_v = \log \sigma_v \xi \quad \mbox{for all} \; \xi \in E_p (k)
\]
we get
\begin{eqnarray*}
  \hat{C} & = & \Big\{ \mu \; \mbox{in} \; \bigoplus_{v \tei \infty} \Z \tei \prod_{v\tei\infty} (\sigma_v \xi)^{\mu_v} \in \mu (\C) \; \mbox{for all} \; \xi \; \mbox{in} \; E_p (k) \Big\} \\
& = & \Big\{ \mu \; \mbox{such that} \; \sum_{v\tei \infty} \mu_v \sigma_v \; \mbox{annihilates} \; E_p (k) \otimes \Q \Big\} \; .
\end{eqnarray*}
Using the isomorphism (\ref{eq:20}) we find
\[
\hat{C} = \Big\{ \mu \; \mbox{such that} \; \sum_{v\tei \infty} \mu_v \sigma_v \; \mbox{annihilates} \; \Q [\Th]^- \Big\} \; .
\]
Set $\varepsilon = \halb (1- c)$. Then $\Q [\Th]^- = \varepsilon \Q [\Th]$ and we see that $\hat{C}$ consists of those $\mu$ such that $\sum_{v \tei \infty} \mu_v (\sigma_v - c\sigma_v)$ annihilates $\Q [\Th]$ i.e. each $\ep \in \Th$. For $\ep , \ep'$ in $\Th$ define the function
\[
\varepsilon_{\ep , \ep'} : \{ v \tei \infty \} \longrightarrow \{ 1 , -1 , 0 \}
\]
by
\[
\varepsilon_{\ep , \ep'} (v) = \left\{ 
  \begin{array}{ll}
1 & \mbox{if} \; \sigma_v \ep = \ep' \\
-1 & \mbox{if} \; \sigma_v \ep = c \ep' \\
0 & \mbox{otherwise} \; .
  \end{array} \right.
\]
Define an inner product on $\bigoplus_{v \tei \infty} \Q$ by the formula
\[
\left\langle \sum a_v \cdot v \; , \; \sum b_v \cdot v \right\rangle = \sum a_v b_v \; .
\]
Now an easy calculation gives the first part of the following theorem if we take into account the relations
\[
\varepsilon_{\ep , \ep'} = \varepsilon_{c\ep , c\ep'} = - \varepsilon_{\ep , c\ep'} = - \varepsilon_{c \ep , \ep'} \; .
\]
We let $S \subset \Th$ be as in proposition \ref{t22}.

\begin{theorem}
  \label{t41}
Let $k / \Q$ be a Galois $CM$ extension. Then we have:
\[
\hat{C} = \left\{ \mu \; \mbox{in} \; \bigoplus_{v \tei \infty} \Z \; \mbox{such that} \; \langle \mu , \varepsilon_{\ep , \ep'} \rangle = 0 \; \mbox{for all} \; \ep , \ep' \in S \right\} \; .
\]
Moreover
\[
\dim \overline{\alpha_{p,\infty} (E_p (k))} = 
\begin{array}{l}
\mbox{dimension of the} \; \Q\mbox{-span of the elements} \; \varepsilon_{\ep , \ep'} \\\
\mbox{for} \; \ep , \ep' \in S \; \mbox{in} \; \bigoplus_{v \tei \infty} \Q \; .
\end{array}
\]
In particular $\alpha_{p,\infty} (E_p (k))$ is dense in $T_k$ if $p$ is completely decomposed.
\end{theorem}

\begin{proof}
  We have
  \begin{eqnarray*}
    \dim \overline{\alpha_{p,\infty} (E_p (k))} & = & \dim \overline{\alpha_{p,\infty} (E_p (k))}^0 = \dim T_k - \dim C \\
& = & \dim T_k - \dim_{\Q} (\hat{C} \otimes \Q) \\
& = & \dim T_k - \dim_{\Q} \langle \varepsilon_{\ep , \ep'} , \ep , \ep' \; \mbox{in} \; S \rangle^{\perp} \\
& = & \dim_{\Q} \langle \varepsilon_{\ep , \ep'} , \ep , \ep' \; \mbox{in} \; S \rangle \; .
  \end{eqnarray*}
If $p$ is completely decomposed, then $\sigma_v \ep = \sigma_w \ep$ implies $v = w$. Moreover $\sigma_v \ep = c \sigma_w \ep$ is not possible for any two places $v, w$. Hence
\[
\varepsilon_{\ep , \ep^{\sigma_w}} (v) = \left\{ 
  \begin{array}{lll}
1 & \mbox{if} & v = w \\
0 & \mbox{if} & v \neq w \; .
  \end{array} \right.
\]
Hence the elements $\varepsilon_{\ep , \ep^{\sigma_w}} \ssp \bigoplus_{v\tei\infty} \Q$ and therefore $\overline{\varepsilon_{p,\infty} (E_p (k))}$ is all of $T_k$.
\end{proof}

\noindent
Mathematisches Institut\\
Westf. Wilhelms-Universit\"at\\
Einsteinstr. 62\\
48149 M\"unster\\
Germany\\
deninge@math.uni-muenster.de
\end{document}